\documentclass[12pt]{article}

\usepackage{amsmath,amssymb,cite}
\font\fmale=cmr8

\def\nn{\nonumber}

\def\bz{{\bar z}}
  
 \def\bv{{\bar v}}
\def\ha{{ {\fmale 1} \over {\fmale 2} }}
\def\om{\omega} 
\def\nt{\noindent}
\font\verysmall=cmr5
\def\phr{\raise1.1pt\hbox{\verysmall x}\kern-8pt\supset}
\def\bhr{\raise1.1pt\hbox{\verysmall x}\kern-9pt\subset}

\def\llra{\longleftrightarrow}

  \def\bbz{Z\!\!\!Z}
\def\bbc{C\kern-6.5pt I}  
\def\a{\alpha}

\def\hp{\hat{\varphi}}  \def\l{\lambda} 

 \def\ce{{\cal E}}

 \def\cy{{\cal Y}}

  \def\ve{\varepsilon}

\def\lra{\longrightarrow}  \def\vr{\vert}

\def\be{{\bar \ce}}   
\def\pd{\partial}   

\def\tB{{\tilde B}}

\def\be{\begin{equation}}
\def\eqn#1{\be\label{#1}}
\def\ee{\end{equation}}

\def\bea{\begin{eqnarray}}
\def\eqnn#1{\bea\label{#1}}
\def\eea{\end{eqnarray}}

\newcommand{\eqna}[1]{\begin{subequations} \label{#1}
\begin{eqnarray}}
\def\eena{\end{eqnarray}
\end{subequations}}

\def\nd{\vfill\end{document}}

\begin{document}

\begin{center}

{\LARGE {\bf Multiparameter Quantum Minkowski Space-Time and Quantum
Maxwell Equations Hierarchy}}

% Classification of Conformal\\[4pt] Representations Induced
%from the \\[6pt] Maximal Cuspidal Parabolic}}

\vspace{10mm}

{\bf \large V.K. Dobrev}

\vspace{5mm}

 Institute for Nuclear
Research and Nuclear Energy \\
{Bulgarian Academy of Sciences}\\ {72 Tsarigradsko Chaussee, 1784
Sofia, Bulgaria}

\end{center}
 \vspace{10mm}

 \begin{abstract}
Earlier we have proposed new $q$ - Maxwell equations which are the
first members of an infinite new hierarchy of $q$ - difference
equations. We have used an indexless formulation in which all
indices are traded for two conjugate variables, $z,\bz$.  We
proposed also new $q$ - Minkowski coordinates which together with
$z,\bz$ can be interpreted as the six local coordinates of a
$SU_q(2,2)$ flag manifold. In the present paper we generalise the
main ingredients of this construction to the multiparameter case
using the seven-parameter quantum group deformation of $GL(4)$ and
$U(gl(4))$ and the four-parameter quantum group deformation of $SL(4)$ and
$U(sl(4))$. The main result is the explicit presentation of the
multiparameter quantum Minkowski space-time within the corresponding
deformed flag manifold.
\end{abstract}

\vspace{10mm}

 \section{Introduction} Invariant differential equations
play a very important role in the description of physical symmetries
- recall, e.g., the examples of Dirac, Maxwell equations, (for more
examples cf., e.g., \cite{BR}).   It is important to construct
systematically  such equations for the setting of quantum groups,
where they are expected as (multiparameter) $q$-difference
equations.

In the present paper we consider the construction of deformed
multiparameter analogs of some conformally invariant equations, in
particular, the Maxwell equations, following the approach of
\cite{Dob}.  We start with the classical situation and we first
write the Maxwell equations in an indexless formulation, trading the
indices for two conjugate variables $z,\bz$. This formulation has
two advantages. First, it is very simple, and in fact, just with the
introduction of an additional parameter, we can describe a whole
infinite hierarchy of equations, which we call the ~{\it Maxwell
hierarchy}~. ~Second, we can easily identify the variables $z,\bz$
and the four Minkowski coordinates with the six local coordinates of
a flag manifold of $SL(4)$ and $SU(2,2)$. Thus, one may look at this
as a nice example of unifying internal and external degrees of
freedom.

Next we need the deformed  analogs of the above constructions. The
specifics of the approach of \cite{Dob} is that one needs also the
complexification of the algebra in consideration. Thus for the $q$ -
conformal algebra we have used the $U_q(sl(m))$ apparatus of
\cite{Dof}  in  the case $m=4$. Thus, in \cite{Domax}  we have
proposed new ~$q$ - {\it Minkowski coordinates}~ as part of the
appropriate $q$ - flag manifold. Using the corresponding
representations and intertwiners of $U_q(sl(4))$ we have also
derived an infinite hierarchy of $q$ - Maxwell equations.

In the present paper we generalise the main ingredients of the above
construction to the multiparameter case. From \cite{DoPar}
 we know that the multiparameter deformation of ~$SL_q(m)$~
and  ~$U_q(sl(m))$~ depends on ~$(m^2-3m+4)/2$~ parameters. We apply
this for ~$m=4$~ in order to consider 4-parameter deformation of the
conformal group, of Minkowski space-time and of Maxwell equations.

\vskip 10mm

\section{Classical setting}

It is well known that Maxwell equations
\eqn{mxa} \pd^\mu  F_{\mu\nu}  ~= ~ J_\nu ~, \quad {\pd^\mu} ^*
F_{\mu\nu}  ~= ~ 0 \ee
or, equivalently
\eqnn{mxb} && \pd_k  E_k
~= ~ J_0 ~(= 4\pi\rho), \quad \pd_0  E_k - \ve_{k\ell m}  \pd_\ell
H_m  ~= ~ J_k ~(=-4\pi j_k), \nn\\  &&\pd_k  H_k  ~= ~ 0 ~, \quad \pd_0
H_k + \ve_{k\ell m}  \pd_\ell  E_m  ~= ~ 0 ~, \eea
where ~$E_k
~\equiv ~F_{k0}$, ~$H_k ~\equiv ~(1/2) \ve_{k\ell m} F_{\ell m}$,
~can be rewritten in the following manner:
\eqn{mxc}  \pd_k  F^\pm_k
~ =~  J_0 ~, \quad \pd_0   F^\pm_k  \pm
 i \ve_{k\ell m}  \pd_\ell  F^\pm_m  ~= ~ J_k ~,  \ee
where
\eqn{mxd} F^\pm_k  ~\equiv  ~E_k   \pm i H_k ~.\ee

Not so well known is the fact that the eight equations in
\eqref{mxc} can be rewritten as two conjugate scalar equations in
the following way: \eqna{mxe} &&I^+ ~ F^+(z)
~=~ J(z,\bz) ~, \\ %&\mxe a\cr
&&I^- ~ F^-(\bz)
~=~ J(z,\bz) ~, \eena %&\mxe b\cr }$$
where
\eqna{ops}
&&I^+ ~=~ \bz \pd_+ + \pd_v  -
\ha\Bigl( \bz z \pd_+ + z\pd_v + \bz \pd_\bv
+ \pd_- \Bigr) \pd_z
~,\\ %&\ops a\cr &&\cr
&&I^- ~=~  z \pd_+ +  \pd_\bv -
\ha \Bigl( \bz z \pd_+ + z\pd_v + \bz \pd_\bv
+ \pd_- \Bigr) \pd_\bz
~, \eena %&\ops b\cr &&\cr
\eqna{nts}
&&x_\pm \equiv x_0 \pm x_3, \quad v \equiv x_1 -i x_2, \quad
\bv \equiv x_1 + i x_2, \\ %&\nts a\cr &&\cr
&& \pd_\pm \equiv \pd/\pd x_\pm, \quad
\pd_v \equiv \pd/\pd v,   \quad \pd_\bv \equiv \pd/\pd \bv,
\eena %&\nts b\cr &&\cr
\eqna{ntsa}
&&F^+(z) ~\equiv ~ z^2 (F_1^+ +iF_2^+) -2z F_3^+ -
(F_1^+ - iF^+_2) ~,\\ %&\ntsa a\cr &&\cr
&&F^-(\bz) ~\equiv ~ \bz^2 (F_1^- -iF_2^-) - 2\bz F_3^- -
(F_1^- +iF^-_2) ~,\\ %&\ntsa b\cr &&\cr
&&J(z,\bz) ~\equiv~ \bz z (J_0 + J_3) + \bz  (J_1 -i J_2) +
z (J_1 +i J_2) +  (J_0 - J_3) ~,\eena % &\ntsa c\cr &&\cr} $$
where we continue to suppress the $x_\mu$, resp., $x_\pm,v,\bv$,
dependence in $F$ and $J$. (The conjugation mentioned above is
standard and in our terms it is : $I^+ \llra I^-$, $F^+(z) \llra
F^-(\bz)$.)

It is easy to recover \eqref{mxc} from \eqref{mxe} - just note that both
sides of each equation are first order polynomials
in each of the two variables $z$ and $\bz$, then
comparing the independent terms in \eqref{mxe} one gets at once
\eqref{mxc}.

Writing the Maxwell equations in the simple form \eqref{mxe} has also
important conceptual meaning. The point is that each of the two
scalar operators ~$I^+,I^-$~
is indeed a single object, namely it
is an intertwiner of the conformal group, while the individual
components in \eqref{mxa} - \eqref{mxc} do not have this interpretation. This is also
the simplest way to see that the Maxwell equations are
conformally invariant, since this is equivalent to the intertwining
property.

Let us be more explicit. The physically relevant representations
$T^\chi$ of the
4-dimensional conformal algebra $su(2,2)$ may be labelled by
$\chi = [n_1,n_2;d]$, where $n_1, n_2$ are non-negative
integers fixing finite-dimensional irreducible representations of
the Lorentz subalgebra, (the dimension being
$(n_1 +1)(n_2 +1)$), and ~$d$~ is the conformal
dimension (or energy). (In the literature these Lorentz
representations are labelled also by $(j_1,j_2)=(n_1/2,n_2/2)$.)
Then the intertwining properties of the
operators in \eqref{ops} are given by:
\eqna{int}
&& I^+ ~:~ C^+ \lra C^0 ~, \quad
I^+ \circ T^+ ~=~  T^0 \circ I^+ ~, \\ %&\int a\cr
&& I^- ~:~ C^- \lra C^0 ~, \quad
I^- \circ T^- ~=~  T^0 \circ I^- ~,\eena  % &\int b\cr}$$
where ~$T^a = T^{\chi^a}$, ~$a=0,+,-$, ~$C^a = C^{\chi^a}$
~are the representation spaces, and the signatures are
given explicitly by:
\eqn{sgn} \chi^+ = [2,0;2] ~,
\quad    \chi^- = [0,2;2] ~, \quad
\chi^0 = [1,1;3] ~,
\ee
as anticipated. Indeed, $(n_1,n_2) = (1,1)$ is the four-dimensional
Lorentz representation, (carried by $J_\mu$ above), and
$(n_1,n_2) = (2,0),(0,2)$ are the two conjugate three-dimensional
Lorentz representations, (carried by $F^\pm_k$
above), while the conformal dimensions are the canonical dimensions
of a current ($d=3$), and of the Maxwell field ($d=2$).
We see that the variables $z,\bz$ are related to the spin
properties and we shall call them 'spin variables'. More
explicitly, a Lorentz spin-tensor $G(z,\bz)$ with signature
$(n_1,n_2)$ is a polynomial in $z,\bz$ of order $n_1,n_2$,
resp. (For more group-theoretical details, cf. \cite{Dob}.)

Formulae \eqref{int}, \eqref{sgn} are part of an infinite hierarchy of
couples of first order intertwiners. Explicitly, instead of
\eqref{int}, \eqref{sgn} we have \cite{Dob}:
\eqna{inta}
&& I^+_n ~:~ C^+_n \lra C^0_n ~, \quad
I^+_n \circ T^+_n ~=~  T^0_n \circ I^+_n ~,\\ %&\inta a\cr
&& I^- _n~:~ C^-_n \lra C^0_n ~, \quad
I^-_n \circ T^-_n ~=~  T^0_n \circ I^-_n ~, \eena %&\inta b\cr}$$
where ~$T^a_n = T^{\chi^a_n}$, ~$C^a_n = C^{\chi^a_n}$, and
the signatures are:
\eqn{sgna} \chi^+_n = [n+2,n;2] ~,
\quad    \chi^-_n = [n,n+2;2] ~, \quad
\chi^0_n = [n+1,n+1;3] ~, \quad  n\in\bbz_+ ~,
\ee
while instead of \eqref{mxe} we have:
\eqna{mxf}
&&I^+_n ~ F^+_n(z,\bz)
~=~ J_n(z,\bz) ~,\\ %&\mxf a\cr
&&I^-_n ~ F^-_n(z,\bz)
~=~ J_n(z,\bz) ~, \eena %&\mxf b\cr }$$
where
\eqna{opsa}
&&I^+_n ~=~ {n+2\over 2} \Bigl(
\bz \pd_+ + \pd_v \Bigr)  -
\ha\Bigl( \bz z \pd_+ + z\pd_v + \bz \pd_\bv
+ \pd_- \Bigr) \pd_z ~, \quad n\in \bbz_+
\qquad\\ %&\opsa a\cr &&\cr
&&I^-_n ~=~  {n+2\over 2} \Bigl( z \pd_+ +  \pd_\bv \Bigr)
- \ha \Bigl( \bz z \pd_+ + z\pd_v + \bz \pd_\bv
+ \pd_- \Bigr) \pd_\bz ~, \quad n\in \bbz_+ \qquad\qquad
\eena %&\opsa b\cr &&\cr}$$
while $F^+_n(z,\bz)$, $F^-_n(z,\bz)$,
$J_n(z,\bz)$, are polynomials in $z,\bz$
of degrees $(n+2,n)$, $(n,n+2)$, $(n+1,n+1)$, resp.,
as explained above.
If we want to use the notation with indices as in
\eqref{mxa}, then $F^+_n(z,\bz)$ and $F^-_n(z,\bz)$ correspond
to $F_{\mu\nu,\a_1,\ldots,\a_n}$ which is
antisymmetric in the indices $\mu,\nu$, symmetric in
$\a_1,\ldots,\a_n$, and traceless in every pair of
indices,\footnote{In 4D conformal field theory the families of
mixed tensors $F_{\mu\nu,\a_1,\ldots,\a_n}$ appear, e.g.,
in the operator product expansion
of two spin $1/2$ fields \cite{DHPS}.}
while $J_n(z,\bz)$ corresponds to $J_{\mu,\a_1,\ldots,\a_n}$
which is symmetric and traceless in every pair of indices.
Note, however, that the analogs of \eqref{mxa} would be much more
complicated if one wants to write explicitly all components.
The crucial advantage of \eqref{mxf} is that the operators ~$I^\pm_n$~
are given just by a slight generalization of ~$I^\pm = I^\pm_0$~.

We call the hierarchy of equations \eqref{mxf}
the ~{\it Maxwell hierarchy}. ~The Maxwell equations are the
zero member of this hierarchy.

Formulae \eqref{mxf},\eqref{inta},\eqref{sgna} are part of a much more general
classification scheme \cite{Dob}, involving
also other intertwining operators, and of arbitrary order. This
scheme was adapted to the $q$-case in  \cite{Dof}.
A  subset of this scheme are two conjugate infinite
two-parameter families of representations which are
intertwined by the same operators \eqref{opsa}.
This  is omitted here for the lack of space, cf. \cite{Dob}.

To proceed further we rewrite \eqref{opsa} in the following
form:
\eqna{opsaz}
&&I^+_n ~=~ \ha \Bigl( (n+2) I_1 I_2  - (n+3) I_2 I_1
\Bigr)
~,\\ %&\opsa a\cr
&&I^-_n ~=~ \ha \Bigl( (n+2) I_3 I_2  - (n+3) I_2 I_3
\Bigr)
~, \eena %&\opsa b\cr}$$
where
\eqn{opsb}
 I_1 ~\equiv ~\pd_z ~, \quad I_2 ~\equiv ~ \bz z \pd_+ + z\pd_v
+ \bz \pd_\bv + \pd_-  ~, \quad I_3 ~\equiv ~\pd_\bz ~. \ee

It is important to note  that group-theoretically the operators ~$I_a$~
correspond to the right action of the three
simple roots of the root system of $sl(4)$, while the operators
~$I^\pm_n$~ are obtained from the lowest possible singular vectors
corresponding to the two non-simple non-highest roots \cite{Dob}.

This is the form that we generalize for the   deformed case.
In fact, we can write at once the general form, which follows
from the analysis of \cite{Dof}:
\eqna{opsc}
&&{\hat I}^+_n ~=~ \ha \Bigl( [n+2]_q {\hat I}_1  {\hat I}_2  -
[n+3]_q {\hat I}_2 {\hat I}_1
\Bigr)
~,\\ %&\opsc a\cr
&&{\hat I}^-_n ~=~ \ha \Bigl( [n+2]_q {\hat I}_3 {\hat I}_2  -
[n+3]_q {\hat I}_2 {\hat I}_3
\Bigr) ~.\eena %&\opsc b\cr}$$
Here ~${\hat I}^\pm_n$~ are obtained from the lowest possible singular vectors
of ~$U_q(sl(4))$, corresponding (as above) to the two non-simple non-highest roots \cite{Dof}.

To proceed further, we should make this form explicit by first generalizing
the variables, then the functions and the operators.

\vskip 5mm

\section{Multiparameter quantum Minkowski space-time}

The variables ~$x_\pm,v,\bv,z,\bz$ ~have definite group-theoretical
meaning, namely, they are six local coordinates on the flag
manifold ~$\cy ~=~ GL(4)/\tB ~=~ SL(4)/B$, where $\tB,B$ are the Borel subgroups of
$GL(4),SL(4)$, respectively, consisting of all upper diagonal matrices.
Under a natural conjugation (cf. also below) this is also a
flag manifold of the conformal group $SU(2,2)$.

We know from \cite{DoPar}
what are the properties of the non-commutative
coordinates on the multiparameter $SL_{q,{\bf q}}$ flag manifold.

There is a technicality here, namely, that we start from the multiparameter deformation
~$GL_{q,{\bf q}}(m)$~  of ~$GL(m)$~ (given by Sudbery \cite{Sud})
which depends on ~$(m^2 -m+2)/2$~ parameters
~$q,q_{ij}$, ~$1\leq i<j\leq m$. (The parametrisation is such that the standard one-parameter deformation is obtained
for all ~$q_{ij}=q$.) Thus, the flag manifold ~$\tilde{\cy}_{q,{\bf q}} ~=~
GL_{q,{\bf q}}(m)/\tB_{q,{\bf q}}(m)$~ depends on the same number of parameters.
For ~$m=4$~ the explicit relations are ($\l \equiv q-q^{-1}$):
\eqnn{coo}
&&x_+ v ~=~ {q_{23}q_{34} \over q_{24}} vx_+ \ ,  \qquad
\bv x_+  ~=~ {q_{14} \over  q_{12}q_{24} }  x_+ \bv\ ,\\ %   &\coo{}\cr
&&x_- v ~=~ {q_{13}\over q_{12}q_{23}  } v x_-\ ,  \qquad
\bv x_- ~=~ {q_{13}q_{34} \over q_{14}} x_-\bv \ ,\nn\\ %  &\cr
&&\bv v ~=~ {q_{13}q_{34} \over q_{12}q_{24}} v\bv \ , \nn\\
&&{q\,q_{24}\over q_{23}q_{34}}x_+ x_- ~=~ {q_{12}q_{24}\over
q\,q_{14}}x_-x_+  ~+~ \l v \bv \ , \nn\eea
\eqnn{coz}
&&\bz z ~=~ {q_{13}q_{24} \over q_{14}q_{23}}z\bz  \  ,\\  %  &\coz{}\cr
&&\bz x_+ ~=~ {q_{13}q_{34} \over q_{14}} x_+\bz \  ,\  \qquad
\bz x_- ~=~ {q_{23}q_{34}\over q^2 q_{24}}x_-\bz  ~+~
\l\bv  \  ,\nn\\
&&\bz \bv  ~=~ {q_{23}q_{34} \over q_{24}} \bv \bz \  ,\  \qquad
\bz v ~=~ {q_{13}q_{34}\over q^2q_{14}}v\bz  ~+~
\l x_+ \  ,\nn\\
&&x_+z ~=~ {q_{14} \over q_{12}q_{24}  } zx_+ \  , \qquad
x_-z ~=~ {q^2 q_{13} \over q_{12}q_{23} }zx_- ~-~
\l v \ , \nn\\
&&v z ~=~ {q_{13}  \over q_{12}q_{23} } z v \  , \qquad
\bv z ~=~ {q^2 q_{14}\over  q_{12}q_{24}}z\bv  ~-~
\l x_+ \ . \nn\eea
Thus, in \eqref{coo} we have a seven-parameter
quantum Minkowski space-time.

We note that when all deformation parameter are phases, i.e., ~$\vr q \vr =1$,
~$\vr q_{ij} \vr =1$, and in addition holds the following relations:
\eqn{relq}
 q_{13} ~=~ {q_{12}q_{24} \over q_{34}} \ , \quad
 q_{14} ~=~ {q_{12}q^2_{24} \over q_{23} q_{34}} \ , \ee
then the commutation relations \eqref{coo} and \eqref{coz} are
preserved by an anti-linear anti-involution
~$\om$~ acting as~:
\eqn{cnj} \om (x_\pm) ~=~ x_\pm ~, \quad \om (v) ~=~ \bv ~, \quad \om (z) ~=~ \bz ~.\ee

Further, we recall from \cite{DoPar} that the dual quantum algebra
~$U_{q,{\bf q}}(gl(m))$~ has the quantum algebra ~$U_{q,{\bf q}}(sl(m))$~ as a commutation subalgebra,
but not as a co-subalgebra. In order to achieve the complete splitting of ~$U_{q,{\bf q}}(sl(m))$~
we have to impose some relations between the parameters, thus the genuine multiparameter deformation
~$U_{q,{\bf q}}(sl(m))$~ depends on ~$(m^2 -3m+4)/2$~ parameters. Using the same conditions we also ensure
that we can restrict from ~$GL_{q,{\bf q}}(m)$~ to ~$SL_{q,{\bf q}}(m)$.

Thus, in the case of ~$m=4$~  for the genuine ~$U_{q,{\bf q}}(sl(4))$~ we have four parameters.
Explicitly, we achieve this by imposing that
the parameters ~$q_{i,i+1}$~ are expressed through the rest as:
\eqn{split}
 q_{12} ~=~ {q^3 \over q_{13}q_{14}} \ ,\quad
q_{23} ~=~ {q^4 \over q_{13}q_{14}q_{24}}\ ,
 \quad q_{34} ~=~ {q^3 \over q_{14}q_{24} }
  \ . \ee %&\split{} }$$

Thus, the four-parameter   quantum Minkowski space-time and
the embedding quantum flag manifold ~$\cy_{q,{\bf q}}$~
are given by \eqref{coo} and \eqref{coz} with
\eqref{split} enforced.

If we would like to enforce also the conjugation \eqref{cnj} then there are more relations between
the deformation parameters, namely, we get:
\eqn{splitz}
q_{12} ~=~ q_{23} ~=~  q_{34} ~=~ {q^2 \over q_{14}}\ , \quad
q_{13} ~=~ q_{24} ~=~  q \ . \ee

Thus, in this case we have a two-parameter deformation and
using the above relations  \eqref{coo} and \eqref{coz} simplify as follows:
\eqnn{coop}
&&x_+ v ~=~ p \, vx_+ \ ,  \qquad
\bv x_+  ~=~ p^{-1} \,  x_+ \bv\ ,\\ %   &\coo{}\cr
&&x_- v ~=~ p^{-1} \,  v x_-\ ,  \qquad
\bv x_- ~=~ p \,  x_-\bv \ ,\nn\\ %  &\cr
&&\bv v ~=~  v\bv \ , \nn\\
&&{q\over p}\, x_+ x_- ~=~ {p \over q }\, x_-x_+  ~+~ \l v \bv \ , \nn\eea
\eqnn{cozz}
&&\bz z ~=~ z\bz  \  ,\\  %  &\coz{}\cr
&&\bz x_+ ~=~ p \, x_+\bz \  ,\  \qquad
\bz x_- ~=~  {p \over q^2 }\, x_-\bz  ~+~
\l\bv  \  ,\nn\\
&&\bz \bv  ~=~ p\, \bv \bz \  ,\  \qquad
\bz v ~=~ {p \over q^2}\, v\bz  ~+~
\l x_+ \  ,\nn\\
&&x_+z ~=~ p^{-1} \, zx_+ \  , \qquad
x_-z ~=~ {q^2 \over p }\, zx_- ~-~
\l v \ , \nn\\
&&v z ~=~ p^{-1} \, z v \  , \qquad
\bv z ~=~ {q^2 \over p }\, z\bv  ~-~
\l x_+ \ , \nn\eea
where ~$p \equiv q^3/q_{14}^2$.

\section{Multiparameter quantum Maxwell equations hierarchy}

The order of variables hinted in \eqref{coo},\eqref{coz} is related to the normal
ordered basis of the quantum  flag manifold ~$\cy_{q,{\bf q}}$~ considered as an
associative algebra:
\eqn{bas} \hp_{ijk\ell mn} ~~=~~
z^i ~v^j ~x_-^k ~x_+^\ell ~\bv^m ~\bz^n ~,
\quad i,j,k,\ell,m,n \in\bbz_+ ~.\ee

\vskip 5mm

We introduce now the representation spaces ~$C^\chi$~, ~$\chi =
[n_1,n_2;d]$~. ~The elements of ~$C^\chi$~, ~which we shall call
(abusing the notion) functions, are polynomials in $z,\bz$ of
degrees $n_1,n_2$, resp., and formal power series in the quantum
Minkowski variables.   Namely, these functions
are given by:
\eqn{ser} \hp_{n_1,n_2}({\bar Y}) ~~=~~ \sum_{{i,j,k,\ell,m,n
\in\bbz_+ \atop i\leq n_1, \ n\leq n_2}} \mu^{n_1,n_2}_{ijk\ell mn}
~\hp_{ijk\ell mn} ~, \ee
where  $\bar Y$ denotes the set of the six coordinates on
$\cy_{q,{\bf q}}$~. ~
Thus the quantum analogs of $F^\pm_n$, $J_n$, cf. \eqref{mxf}, are~:
\eqn{fnc} {\hat F}^+_n ~=~ \hp_{n+2,n}({\bar Y}) ~, \quad
{\hat F}^-_n ~=~ \hp_{n,n+2}({\bar Y}) ~, \quad
{\hat J}_n ~=~ \hp_{n+1,n+1}({\bar Y}) ~.\ee

\vskip 5mm

Using the above machinery we can present a deformed  version
of the Maxwell hierarchy of equations. First, we mention that the
explicit form of the operators  ~$I_a$~ in \eqref{opsb} is obtained by
the infinitesimal right action of the three simple root
generators of $sl(4)$ on the flag manifold $\cy$ (following the
procedure of \cite{Dob}). In the deformed case the right action of
$U_{q,{\bf q}}(sl(4))$ on $\cy_{q,{\bf q}}$ is known from \cite{DoPar}, thus,
we have:
\eqn{larcq} {\hat I}_a ~=~ \pi_R (X^-_a) \ee

From this we obtain the multi-parameter quantum   Maxwell hierarchy of equations by
  substituting the operators of \eqref{larcq} in \eqref{opsc}, i.e.,
  the final result is:
\eqna{mxg}
&&{\hat I}^+_n ~ {\hat F}^+_n
~=~ {\hat J}_n ~,\\ %&\mxg a\cr
&&{\hat I}^-_n ~ {\hat F}^-_n
~=~ {\hat J}_n ~. \eena % &\mxg b\cr }$$
The reason that we can use \eqref{opsc} is that the multiparameter $U_{q,{\bf q}}(sl(4))$
depends only on ~$q$~ as a commutation subalgebra, while the dependence on the other parameters
is exhibited only in its co-algebra structure and in the explicit expressions of
~$\pi_R (X^-_a)$.

\vskip 5mm

\nt {\bf Remark:} Certainly, as we did in the one-parameter case \cite{Domax}, we would like
to present \eqref{larcq} and \eqref{mxg} more explicitly, cf. \cite{Doprep}.

\vskip 10mm

\section*{Acknowledgments} The author has received partial support from COST Actions
MP1210 and MP1405, and from Bulgarian NSF Grant DFNI T02/6.

\end{document}